\documentclass{article}
\usepackage{graphicx}
\usepackage{amsmath}
\usepackage{amssymb}
\usepackage{derivative}
\usepackage{braket}
\usepackage{comment}
\usepackage{url}
\usepackage{subcaption}
\usepackage{caption}
\usepackage{svg}
\usepackage{hyperref}
\usepackage{authblk}

\title{Coupled Arnol'd cat maps on circulant graphs}

\author[1]{K. Manolas \thanks{Corresponding author. E-mail: kmanolas@phys.uoa.gr}}

\author[1,2]{E. Floratos}

\affil[1]{Faculty of Physics, Department of Nuclear and Particle Physics, National and Kapodistrian University of Athens, University Campus, Ilisia, 15771 Athens, Greece}
\affil[2]{Academy of Athens, Panepistimiou 28, 10679 Athens, Greece}
\date{}

\begin{document}

\maketitle

\abstract{This paper investigates the chaotic properties of Arnol'd cat maps (ACMs) coupled on the nodes of a circulant graph. By demanding that the system's evolution matrix be symplectic, we determine the coupling matrix, which is naturally interpreted as the adjacency matrix of a circulant graph. Specifically, the study analyses the system's Lyapunov spectra and Kolmogorov-Sinai (K-S) entropy. Numerical simulations yield the counterintuitive result that the entropy production does not increase as the connectivity of the graph increases, due to the translational symmetry of the circulant graph. Moreover, we analyse the spectra of the periods of the evolution matrix on a finite toroidal phase space of the dynamical system.

\vspace{1em}
}
\newpage

\tableofcontents

\section{Introduction}
Graphs serve as powerful and versatile structures for modeling and understanding dynamical systems, with applications in areas such as those of machine learning, sociology, finance, neuroscience, social networks, and physical sciences.

A graph is a combinatorial object often defined as $\mathcal{G}=(V,E)$ where $V$ a set of nodes (vertices) and $E$ the set of edges (links) connecting the nodes. Assuming some ordering of the nodes, a graph's connectivity can be represented through a square matrix $A\in \mathbb{R}^{|\mathcal{V}| \times |\mathcal{V}|}$, called an \textit{adjacency matrix}. There are two important graph classifications based on the adjacency matrix $A$; unweighted ($Aij \in {0,1})$ and weighted graphs, and undirected (where $A$ is symmetric) and directed graphs.
A specific subset of graphs, with great theoretical and practical interest, is that of \textit{circulant graphs}, whose connectivity is described by a circulant matrix\footnote{More precisely, a graph is called circulant if there exist any permutation of its nodes for which the adjacency matrix is circulant. In our case, as we consider the adjacency matrix to be the more fundamental object and construct a graph through it, we forgo any mention to permutations of the nodes and simply call the graph circulant if the adjacency matrix we use for its construction is circulant.}, that is, a matrix whose lines (rows) are constructed through consecutive circular shifts of a generating vector (see for example \cite{Gray2001}). Consequently, every node of a circulant graph has the same connectivity making the graph rotation invariant for any rotation that is an integer multiple of $2 \pi / |V|$. 

Circulant matrices have the useful property of being diagonalisable by the Discrete Fourier Transform, allowing both for quick computations by moving back and forth between the configuration and frequency domains, as well as providing one of the few topologies that are analytically solvable. Furthermore, the action of a circulant matrix on a signal vector corresponds to a circular convolution with a finite filter (again, see \cite{Gray2001}), which has grown to be of special interest with the tremendous success of Convolutional Neural Networks in Machine Learning \cite{goodfellow2016}.

As a result of the analytic properties of their adjacency matrices, circulant graphs provide a type of toy model bridging the gap between theory and practice, where a number of standard signal processing operations can be defined in a natural way \cite{Shuman2013}. Such graphs have seen increased use in the field of machine learning with the rise of Geometric Deep Learning \cite{bronstein2021geometric}, where their inherent rotational invariance is leveraged \cite{Huang2023}, as well as a means to compress the weight matrices' sizes and expedite their training process through the use of circulant, or block circulant, weight matrices \cite{Cheng2015, Ding2017}. There have also been applications of circulant graphs in the context of quantum mechanics, where it has been shown that integral circulant graphs allow for perfect quantum state transfer \cite{BASIC2009}.

Another field where graphs arise naturally, and have been employed extensively, is that of dynamical systems. By placing maps on the nodes of a graph and coupling them together, Coupled-Map Lattices (CMLs) \cite{Kaneko1992} provide a low-computational-cost class of models that allow for the study of a wide variety of chaotic systems. They were first introduced as a model for spatio-temporal chaos \cite{Kaneko1984}, and were shown to exhibit a plethora of characteristics, especially after their generalisation from simple grids to graphs, such as synchronisation (see for example \cite{Gade2000}). Circulant, or block circulant, graphs were one of the first topologies to be studied \cite{Amritkar1991}, again due to the ease of diagonalising their adjacency matrix through the DFT. The literature is mainly concentrated around ring graphs \cite{Waller1984}, the simplest class of circulant graphs, were the theoretically predicted chimera states \cite{Omelchenko2011} have been observed in practice \cite{Hagerstrom2012}.

Although CMLs have historically made use of chaotic maps, there have been some attempts to better understand them by employing linear maps, as they are analytically tractable. Since a purely linear map cannot be chaotic, chaos is usually introduced either through piecewise linear maps (e.g. \cite{Just1998}), or by imposing periodic boundary conditions on the lattice (see e.g. \cite{Gutkin_2016, Axenides2023}). Understanding the dynamical properties of such chaotic systems is of growing practical importance. As such systems are analytically solvable, being linear, they provide highly tractable models that could serve as efficient computational reservoirs for training machine learning models \cite{TANAKA2019}. Furthermore, the mixing properties resulting from their chaotic dynamics suggests potential future applications in Graph Neural Networks, to avoid optimisation issues such as vanishing gradients \cite{Tao2026}.

In reference \cite{Axenides2023}, the authors constructed another linear CML model on the torus, employing the Arnol'd cat map (ACM) and using its connection to the Fibonacci sequence to allow for analytic calculations. By demanding the evolution matrix be symplectic, they showed that the connectivity matrix must be symmetric, and provided the general form of its blocks, assuming translation invariance.  Building upon that work, we demonstrate that mapping these systems onto circulant graphs, reveals the novel physical phenomenon of highly symmetric connectivities inducing destructive interference which supresses Kolmogorov-Sinai entropy production.

The plan of the paper is as follows:

In Section~\ref{fib_main_section} we provide a brief overview of the Arnol'd cat map and its connection to the Fibonacci sequence.

In Section~\ref{section_coupled_acm} we generalise to any number of coupled ACMs following the construction of \cite{Axenides2023}.

In Section~\ref{Sect_chaos_circ_graph} we specialise the circulant coupling matrix to be the adjacency matrix of a circulant graph. We then study the chaotic properties of the coupled ACMs on this graph by analytic investigation of the Lyapunov spectra, Kolmogorov-Sinai entropy, and period spectra of the system.

In Section~\ref{section_conclusion} we summarise our results and discuss open issues and potential applications.

\section{Fibonacci sequence and the Arnol'd cat map}\label{fib_main_section}
We begin with a quick overview of the Arnol'd cat map and its action on the phase space of a single particle, as well as its connection to the Fibonacci sequence.
We choose the phase space to be a periodic lattice of lattice of length equal to 1, and the Arnol'd cat map
\begin{equation}\label{acm_def}
    A = \begin{bmatrix}
        1 & 1 \\
        1 & 2
    \end{bmatrix}
\end{equation}
acts as
\begin{equation}\label{acm_action_points}
    (r_m, s_m) A = (r_{m+1}, s_{m+1}) \implies (r_{m+1}, s_{m+1}) = (r_m + s_m, r_m + 2s_m) \ mod1
\end{equation}
with the calculations being $mod1$, so that the coordinates of the lattice points
\begin{equation}\label{amc_lattice_point_evol}
    \begin{split}
        r_{m+1} &= r_m + s_m\ mod1 \\
        s_{m+1} &= r_m + 2s_m\ mod1
    \end{split}
\end{equation}
take values in $[0,1]$.

If $r_0, s_0$ are rational they can always be written as $r_0 = k/N, \ s_0 = l/N \leq 1$ leading to
$r_1 = \frac{k+ l}{N}$, $s_1 = \frac{k+2l}{N}$, were calculations are again $mod1$. We can then rewrite the movement on the lattice as
\begin{align}
    r_{m+1} &= k_m + l_m \ modN \\
    s_{m+1} &= k_m +2l_m \ modN
\end{align}
with $k,l = 0, 1, \dots, N-1$, representing the movement on the double toroidal space $\mathbb{Z}_N \times \mathbb{Z}_N = \mathbb{T}_{2N}$.

Since $(r_{m+1}, s_{m+1}) = (r_m, s_m)A \ modN$ it is apparent that $(r_m, s_m) = (r_0, s_0)A^m$ $modN$ for $m =0, 1, \dots$ . Given that the lattice is finite, there must be a period $T$ for the map $A$, such that $A^T = I_{2 \times 2} \ modN$. This map depends on $N$ in a random way, as shown in \cite{Dyson1992}.

In order to make the relation of the Arnol'd cat map and the Fibonacci sequence apparent, let us first present a short overview of well known results surrounding the Fibonacci sequence. We denote as $f_m$ the Fibonacci integers satisfying the sequence and initial conditions
\begin{equation}\label{fib_seq}
    f_{m+1} = f_m + f_{m-1} , \quad f_0 = 0,\ f_1 = 1
\end{equation}
which can be written in matrix form as
\begin{equation}\label{fib_seq_matr}
    (f_m , f_{m+1}) = (f_{m-1}, f_m) 
    \begin{pmatrix}
        0 & 1 \\
        1 & 1
    \end{pmatrix}
\end{equation}
with $\begin{pmatrix}
        0 & 1 \\
        1 & 1
    \end{pmatrix} \begin{pmatrix}
        0 & 1 \\
        1 & 1
    \end{pmatrix} = \begin{pmatrix}
        1 & 1 \\
        1 & 2
    \end{pmatrix} = A$ the Arnol'd cat map for $n=1$. It is important to note that the matrix of the Fibonacci sequence is \textit{anti-symplectic}, that is it satisfies
    \begin{equation}\label{fib_matr_anti_symp}
        \begin{pmatrix}
        0 & 1 \\
        1 & 1
    \end{pmatrix}^T \epsilon \begin{pmatrix}
        0 & 1 \\
        1 & 1
    \end{pmatrix} = - \epsilon
    \end{equation}
    with $\epsilon = \begin{pmatrix}
        0 & -1 \\
        1 & 0
    \end{pmatrix}$ the symplectic form, and as such is not symplectic, but its square is an element of the symplectic group $SL_2(\mathbb{Z})$, that is
    \begin{equation}\label{acm_symp}
        A^T \epsilon A =\epsilon \ .
    \end{equation}
    
    The Fibonacci sequence can also be written as
    \begin{equation}\label{fib_seq_power}
        (f_m , f_{m+1}) = (0, 1) 
    \begin{pmatrix}
        0 & 1 \\
        1 & 1
    \end{pmatrix}^m = \begin{pmatrix}
        0 \\
        1
    \end{pmatrix}
    \begin{pmatrix}
        f_{m-1} & f_m \\
        f_m & f_{m+1}
    \end{pmatrix}
    \end{equation}
    for $m=1,2,\dots$ . Consequently the time evolution under the ACM can be written as:
    \begin{equation}\label{acm_power_form}
        A^m = \begin{pmatrix}
            0 & 1 \\
            1 & 1
        \end{pmatrix}^{2m}
        = \begin{pmatrix}
            f_{2m-1} & f_{2m} \\
            f_{2m} & f_{2m+1}
        \end{pmatrix}
    \end{equation}
We can solve the recursion relation for the Fibonacci sequence \eqref{fib_seq} by using the ansatz  $f_m = a \rho^m$. Solving the resulting equation $\rho^2 - \rho -1 = 0$ we get the well known solutions
\begin{equation}\label{fib_eigenva}
    \rho_{\pm} = \frac{1 \pm \sqrt{5}}{2} \ .
\end{equation}
Given that the recursion relation is linear with coefficients independent of $m$, the linear combination $f_m = a\rho_+^m + \beta \rho_-^m$ is also a solution. By making use of the initial conditions $f_0=0$ and $f_1 = 1$ we find:
\begin{equation}\label{fib_int_solved}
    f_m = \frac{\rho_+^m - \rho_-^m}{\sqrt{5}}
\end{equation}

As a result, the movement on the lattice is described by
\begin{equation}\label{acm_power_lattice_act}
    (r_m, s_m) = (r_0, s_0) \begin{pmatrix}
            f_{2m-1} & f_{2m} \\
            f_{2m} & f_{2m+1}
        \end{pmatrix} \xrightarrow[m \rightarrow \infty]{} f_{2m}
        \begin{pmatrix}
            \rho_+ & 0 \\
            0 & \frac{1}{\rho_+}
        \end{pmatrix}
\end{equation}
where it is noted that $\rho_{\pm}$ are the eigenvalues of the Fibonacci matrix and their squares are the eigenvalues of the Arnol'd cat map (ACM). Thus the ACM expands the phase space along the eigenvector corresponding to $\rho_+$ and contracts the phase space along the eigenvector corresponding to $\rho_-$, leading to strong chaos with the positive Lyapunov exponent \cite{arnold1989} equal to $l = \ln{\rho_+}$.

\section{Coupled Arnol'd cat map}\label{section_coupled_acm}
In this section we generalise the above construction to any number of interacting ACMs. We begin with two particles (two ACMs), using two interacting Fibonacci sequences with added interacting terms
\begin{equation}
\begin{split}\label{double_acm_int}
    f_{m+1} &= a_1 f_m +b_1 f_{m-1} + c_1 g_{m} + d_1 g_{m-1} \\
    g_{m+1} &= a_2 g_m +b_2 g_{m-1} + c_2 f_{m} + d_2 f_{m-1}
\end{split}
\end{equation}
where $a_i, b_i, c_i, d_i \in \mathbb{Z}$ for $i = 1,2$. The above can be written in matrix form as
\begin{equation}\label{block1}
X_{m+1}= \begin{pmatrix}
    f_m \\
    g_m \\
    f_{m+1} \\
    g_{m+1}
\end{pmatrix} =
\left(\begin{array}{cc} 0_{2\times 2} & I_{2\times 2} \\ D & C\end{array}\right) \begin{pmatrix}
    f_{m-1} \\
    g_{m-1} \\
    f_m \\
    g_m
\end{pmatrix} = L X_{m}
\end{equation}
with 
\begin{equation}\label{block_compon}
    D = \begin{pmatrix}
        b_1 & d_1 \\
        d_2 & b_2
    \end{pmatrix}, \quad
    C = \begin{pmatrix}
        a_1 & c_1 \\
        c_2 & a_2
    \end{pmatrix}
\end{equation}
Similarly to the case of a single ACM, we demand the Fibonacci matrix (in this case $L$ \eqref{block1}) be anti-symplectic
\begin{equation}\label{fib_n_anti_sym}
    L J L = -J
\end{equation}
where $J = \begin{pmatrix}
    0_{n\times n} & -I_{n\times n}\\
    I_{n\times n} & 0_{n\times n}
\end{pmatrix}$ is the symplectic form in the 2n-dimensional toroidal phase space of unit length (the 2n-dimensional hypercube with periodic boundary conditions), so that its square, the ACM for any $n$, is symplectic.

The constraint \eqref{fib_n_anti_sym} implies 
\begin{equation}\label{fib_constr}
    D = I_{n \times n} \ , \quad C = C^T
\end{equation}
resulting in $b_1 = b_2 = 1$, $d_1 = d_2 = 0$, $c_1 = c_2 \equiv c$ and $a_1 = k_1$, $a_2 = k_2$, that lead to the Fibonacci matrix
\begin{equation}\label{Fib_matrix}
    L = \begin{pmatrix}
        0 & 0 & 1 & 0 \\
        0 & 0 & 0 & 1 \\
        1 & 0 & k_1 & c \\
        0 & 1 & c & k_2 
    \end{pmatrix}
\end{equation}
corresponding, through \eqref{block1}, to the coupled Fibonacci sequences
\begin{equation}
\begin{split}\label{double_acm_symp}
    f_{m+1} &= k_1 f_m + f_{m-1} + c g_{m} \\
    g_{m+1} &= k_2 g_m + g_{m-1} + c f_{m}
\end{split}
\end{equation}

It is now apparent, that for the general case of $n$ coupled Fibonacci sequences, the matrix is of the form
\begin{equation}\label{fib_matr_n}
    L = \begin{pmatrix}
        0_{n \times n} & I_{n \times n} \\
        I_{n \times n} & C
    \end{pmatrix}
\end{equation}
again satisfying $L^T J L = -J$, so that its square
\begin{equation}\label{symp_matrix_n}
    M = L^2 = \begin{pmatrix}
        I_{n \times n} & C \\
        C & I_{n \times n} + C^2
    \end{pmatrix}
\end{equation}
belongs in the group $Sp_{2n}[\mathbb{Z}]$, that is it satisfies $M^T J M = J$, and as a result its eigenvalues come in pairs $(\lambda, 1/\lambda)$, with $\lambda > 1$. Defining the diagonal integer matrices $G_0$ and $G_1$, as well as the circulant matrix $P = \delta_{i, j+1}$, the coupling matrix for $n$ sequences, each only coupled to their nearest neighbour, can be written as
\begin{equation}\label{coupl_matr_non_trans}
    C = G_0 + PG_1 + G_1P^T
\end{equation}
where the orthogonal matrix $P$ can be identified as the translation operator on the double toroidal lattice $T_N = \mathbb{Z}_N \times \mathbb{Z}_N$ in Finite Quantum Mechanics \cite{ATHANASIU1994}. We are interested in studying the translation invariant case of the coupling, that is to say, for diagonal matrices that are multiples of the identity matrix $G_0 = g_0 \delta_{i,j}$ and $G_1 = g_1 \delta_{i,j}$,
\begin{equation}\label{near_neighb_coupl_matr}
    C = G_0 + G_1(P + P^T)
\end{equation}
which in turn is equivalent to the coupling matrix $C$ commuting with the matrix $P$. Since $C$ is symmetric, so is $L$, and consequently $M = L^2$ is positive definite.

The coupling matrix can be generalised to include more than just nearest neighbour interactions, simply by including powers of the matrix $P$,
\begin{equation}\label{gen_transl_coupl_matr}
    \begin{split}
    C &= g_0 I_{n \times n} + \sum_{l=1}^{n-1} g_l P^l \\ &= g_0 I_{n \times n} +\sum_{l=1}^{[n/2]}g_l P^l + \sum_{l=[n/2]+1}^{n-1}g_l P^l \\ &= g_0 I_{n \times n} + \sum_{l=1}^{[n/2]} g_l \left(P^l + \left(P^T\right)^l\right)
    \end{split}
\end{equation}
where the term $g_i \left(P^i + \left(P^T\right)^i\right)$ corresponds to the $i$-th neighbour interaction, and the equation $P^i = (P^T)^{n-i}$ must hold so that the coupling matrix $C$ remains symmetric. The block form of the evolution matrix remains the same
\begin{equation}\label{evol_matrix_general}
    M = L^2 = \begin{pmatrix}
        I_{n \times n} & C \\
        C & I_{n \times n} + C^2
    \end{pmatrix}
\end{equation}
with $C$ now being the matrix for $n$, fully coupled, ACMs.

The dynamical system of n linearly interacting ACMs is described by the $2n \times 2n$ symplectic evolution matrix $M$, acting on the toroidal phase space of length $1$ (all matrix operations are modulo 1), that is 
\begin{equation}\label{system_evol}
    x_{m+1} = x_m M
\end{equation}
where $x_m = (q_m, p_m)$, and $q_m, p_m$ are the positions and momenta of the n particles at time step $m$. Because $M$ is a symplectic matrix it preserves the volume of the phase space.

Given that the coupling matrix \eqref{near_neighb_coupl_matr}, \eqref{gen_transl_coupl_matr} is circulant, being a sum of circulant matrices, it is diagonalised by the Finite Fourier Transform
\begin{equation}\label{FFT}
    F^\dagger C F = D \quad \text{where} \ F_{r,s} = \frac{1}{\sqrt{n}} e^{2 \pi i rs/n} \equiv \frac{1}{\sqrt{n}} \omega^{rs}_n
\end{equation}
with its eigenvectors being the columns of the FFT $v_j = (1, \omega^j_n, \omega^{2j}_n, \dots, \omega^{(n-1)j}_n)^T$ and the corresponding eigenvalues being determined by the first row of the matrix $C$, i.e. its generating vector
\begin{equation}\label{gen_vector}
    \begin{split}
        C_{j,0} &= (g_0, g_1, g_2, \dots, g_{(n-1)/2}, g_{(n-1)/2}, \dots, g_2, g_1) \quad \text{for}\ n= \text{odd} \\
        C_{j,0} &= (g_0, g_1, g_2, \dots, g_{n/2-1}, g_{n/2}, g_{n/2-1}, \dots, g_2, g_1) \quad \text{for}\ n= \text{even}
    \end{split}
\end{equation}
as
\begin{equation}\label{coupl_eigenva_odd}
    d_j = g_0 + \sum_{l=1}^{(n-1)/2} g_l \left( e^{\frac{2 \pi i j l}{n}} + e^{\frac{2 \pi i j (n-l)}{n}} \right) = g_0 + 2 \sum_{l=1}^{(n-1)/2} g_l \left( \cos{\frac{2 \pi j l}{n}} \right)
\end{equation}
for $n =$ odd, since $\exp\left(\frac{2 \pi i j (n-l)}{n}\right) = \exp\left(-\frac{2 \pi i j l}{n}\right)$, and
\begin{equation}\label{coupl_eigenva_even}
    d_j = g_0 + 2 \sum_{l=1}^{n/2 -1} g_l \left( \cos{\frac{2 \pi j l}{n}} \right) + g_{n/2}(-1)^j
\end{equation}
for $n=$ even.

\section{Arnol'd cat maps on a circulant graph}\label{Sect_chaos_circ_graph}
In this section we shall connect the matrix $C$, which is circulant, with the adjacency matrix of a circulant graph, and we will proceed to the study of the chaotic properties of the system of $n$ coupled ACMs on the nodes of such graph. We remind that a circulant graph is a graph with fixed connectivity, same for every node. Moreover it has a built-in rotational invariant structure which is appropriate for embedding symmetries on graphs. We shall use only unweighted graphs in this section although our discussion is easily extendable to weighted graphs as well. The coupling vector of \eqref{evol_matrix_general} $(g_0, g_1, \dots, g_{n-1})$ will be a binary vector expressing simply the connectivity between the nodes, so the components take values in $\{0,1\}$. Furthermore, to study symplectically interacting nodes, the adjacency matrix must be symmetric, and as such, a binary vector will be of the form $g = (g_0, g_1, g_2,$$\dots, g_{n/2-1},$$ g_{n/2}, g_{n/2-1}, $$\dots, g_2, g_1)$\footnote{For simplicity's sake we have assumed $n$ to be even.}.

With the connection between $n$ coupled ACMs and an $n$-node circulant graph evident, we can now proceed with the study of a graph's chaotic properties by studying those of the corresponding system of coupled ACMs. To do so, we will consider the system's Lyapunov spectra and Kolmogorov-Sinai entropy, since the eigenvalues are easily obtainable through the DFT. Henceforth, we will use the circulant matrix's generating vector $g = (g_0, g_1, \dots, g_{n-1})$, $g_i \in \{0,1\}$, to describe a circulant graph and its respective system of coupled ACMs.

Having constructed the evolution matrix $M$ of the system, we must describe the phase space of the system on which this matrix acts. On every node of the graph we assign a 2-dimensional toroidal phase space resulting in the phase space of the system being the $2n$-dim torus. The points $x$ of this phase space describe the positions and momenta of the $n$ particles, each one assigned to a node of the graph. The classical discrete equations of motions can be written using the evolution matrix $M$ from \eqref{amc_lattice_point_evol} as
\begin{equation}\label{class_eoms}
    x_{m+1} = x_m M \implies (q_{m+1}, p_{m+1}) = (q_m, p_m) \begin{pmatrix}
        I_{n \times n} & C \\
        C & I_{n \times n} + C^2
    \end{pmatrix}
\end{equation}
where $m=0,1,2,\dots$ the timestep of the map, and $q_m$ and $p_m$ the positions and momenta of $n$ particles. Specifically, $q_m, p_m$ are $n$-dimensional vectors whose components take values between 0 and 1 and the evolution equation \eqref{class_eoms} is computed $mod1$. Next we solve the equations of motion
\begin{align}
    q_{m+1} &= q_m + p_m C \label{pos_eom} \\
    p_{m+1} &= q_m C + p_m ( I_{n \times n} + C^2) \label{mom_eom}
\end{align}
in terms of positions only, by solving \eqref{pos_eom} in terms of momentum and substituting in \eqref{mom_eom}, giving Newton's discrete equations of motion
\begin{equation}\label{pos_only_eom}
    q_{m+1} - 2q_m + q_{m-1} = q_m C^2
\end{equation}
in the case where $C^{-1}$ exists, though otherwise the same result can still be obtained, albeit with a few more algebraic steps.

As it was shown in Section \ref{section_coupled_acm}, the matrix $C$ is diagonalised by the FFT. Consequently by applying the FFT on the equations of motion \eqref{pos_only_eom}, we get
\begin{equation}\label{diag_pos_eom}
    r_{m+1} - 2r_m + r_{m-1} = r_m D^2
\end{equation}
where we have defined $r_m \equiv q_m F$, with $F$ being the FFT and $D$ is the diagonal form of the matrix $C$. Since the matrix $C^2$ by its construction is a square integer matrix, its eigenvalues are non-negative. Similarly to the approach in Section \ref{fib_main_section}, in order to solve \eqref{diag_pos_eom}, we make the ansatz $(r_m)_i = \delta_{i,j} \rho_i^m a_j$, $i,j=0,1,\dots,n-1$, for the $i$-th component of the of the normal mode, at time step $m=0,1,2,\dots$, resulting in
\begin{equation}\label{normal_modes_rho}
    \rho_i^2 - (2 + d_i^2) \rho_i + 1 = 0 \implies \rho_{\pm, i} = \frac{2 + d_i^2}{2} \pm \frac{|d_i|}{2} \sqrt{d_i^2 + 4}
\end{equation}
which in turn allows for the general solution of $r_m$ to be written as:
\begin{equation}\label{normal_modes_sol}
    r_m = a_+ \rho_+^m + a_- \rho_-^m
\end{equation}
The constant vector can be determined by the use of initial conditions $(q_0, q_1)$, and the normal modes using the equation $q_m = r_m F^\dagger$. 

The Lyapunov exponents of a dynamical system, $\lambda_k$, $k=0,1,\dots,n-1$, characterize the rate of separation of nearby trajectories in phase space. For $n$ coupled ACMs on the nodes of a symmetric circulant graph, the positive Lyapunov exponents are given by the solution $\rho_{+}$, of \eqref{normal_modes_rho}, as
\begin{equation}\label{lyap_exp}
    \lambda_{+,k} = \ln{\rho_{+}} = \ln{\left(\frac{2 + d_k^2}{2} + \frac{|d_k|}{2} \sqrt{d_k^2 + 4}\right)}, \quad k=0,1,\dots, n-1
\end{equation}
with the eigenvalues $d_k$ of the matrix $C$ given by \eqref{coupl_eigenva_odd}, \eqref{coupl_eigenva_even} for odd and even $n$ respectively.

Building on this, the Kolmogorov-Sinai entropy \cite{Pesin1977, Sinai1989}, which corresponds to the rate of entropy production, is the sum of all positive Lyapunov exponents
\begin{equation}\label{ks_entropy}
    S_{KS} = \sum_{i=0}^{n-1} \lambda_{+,i}
\end{equation}
which must be a linear function of $n$ for large $n$.

\begin{figure}[htbp]
  \centering
  \begin{subfigure}[b]{0.5\textwidth}
    \centering
    \includegraphics[width=\textwidth]{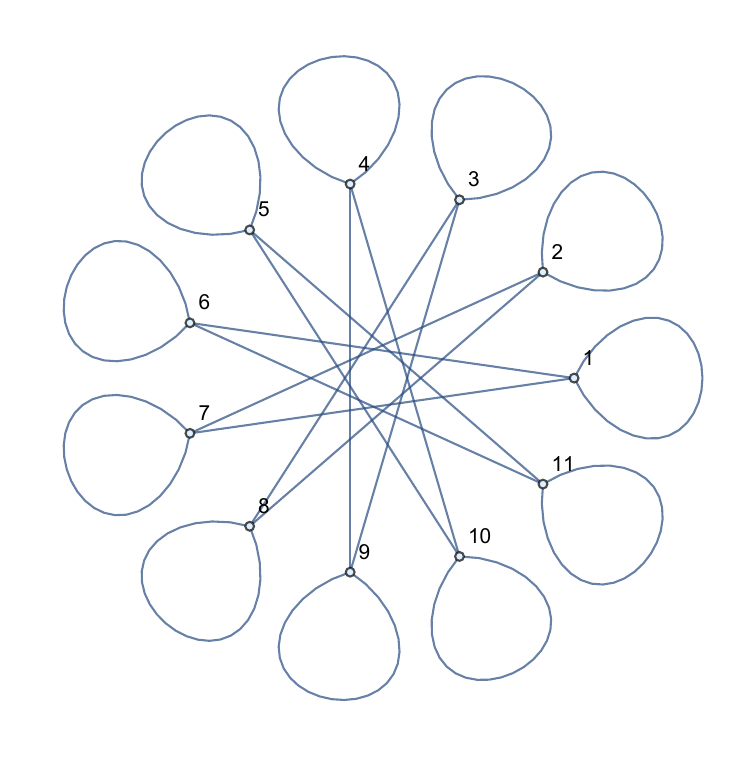}
    \caption{Circulant graph dictated by the adjacency matrix generated by $g$.}
    \label{fig1_subfig4}
  \end{subfigure}
  \hfill
  \begin{subfigure}[b]{0.475\textwidth}
    \centering
    \includegraphics[width=\textwidth]{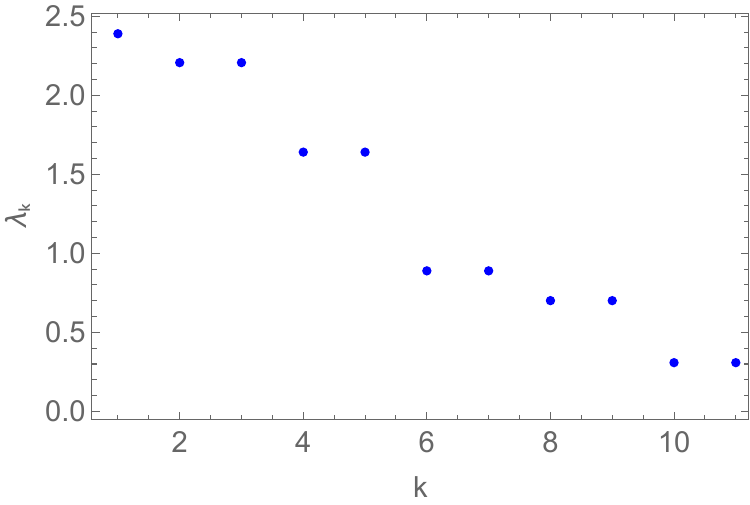}
    \caption{Sorted positive Lyapunov exponents of $g$ by descending magnitude.}
    \label{fig1_subfig2}
  \end{subfigure}

  \caption{In (a) the graph for the system encoding the integer 1072, that is for the adjacency matrix generated by $g=(1,0,0,0,0,1,1,0,0,0,0)$, is shown as an example, noting the rotational invariance characterising all such binary circulant graphs. The positive Lyapunov exponents are  shown in (b), having been sorted in order of descending magnitude.}
  \label{figure_1}
\end{figure}
In order to give some concrete numerical results we shall first look at two systems, one described by a vector of size $n=11$, small enough that the graph can be feasibly visualised with a plot, and another described by a vector of increasing length from $n=3,5,7,\dots,101$, to visualise how the system scales for larger $n$, through the plot of the K-S entropy.

It is known that a positive integer $m$, can be decomposed into its unique binary form
\begin{equation}\label{binary_form_int}
    m = x_k 2^k + x_{k-1} 2^{k-1} + \dots + x_1 2 + x_0
\end{equation}
where $x_i \in \{0,1\} \ \forall i$, and $k$ such that $2^k \leq m < 2^{k+1}$. It is clear from the above decomposition that the largest bit $x_k$ will always equal 1, and the last bit $x_0$ will equal 1 for $m=odd$ and 0 for $m=even$\footnote{Clearly the decomposition could be made non-unique by keeping the coefficients $x_l$, $l>k$, which would all equal zero.}. As a specific example, we take the integer $m=1072$, which decomposes into $m=2^{11} + 2^6 + 2^5$. The decomposition implies that the bits $x_{11} = x_6 = x_5 = 1$ and the rest equal zero, or in vector form $x = (1,0,0,0,0,1,1,0,0,0,0)$. Consequently, by taking $x$ to be the generating vector of a circulant adjacency matrix with matrix elements $A_{i,j} = g_{(j-i)modn}$, we can create a circulant graph with the integer $m$ encoded in its structure. The plots of the graph of the system, and the corresponding sorted positive Lyapunov exponents are shown in Figure~\ref{figure_1}.

\begin{figure}[htbp]
  \centering
  \begin{subfigure}[b]{0.49\textwidth}
    \centering
    \includegraphics[width=\textwidth]{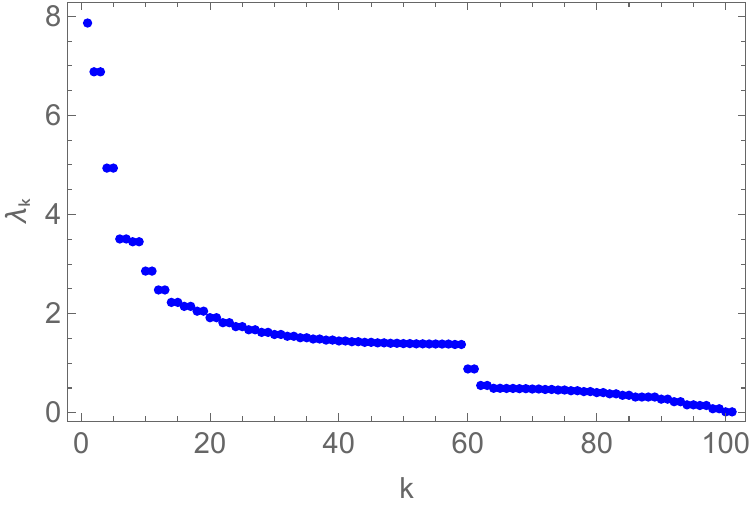}
    \caption{Sorted positive Lyapunov exponents of $h$ by descending magnitude.}
    \label{fig2_subfig2}
  \end{subfigure}
  \hfill
  \begin{subfigure}[b]{0.5\textwidth}
    \centering
    \includegraphics[width=\textwidth]{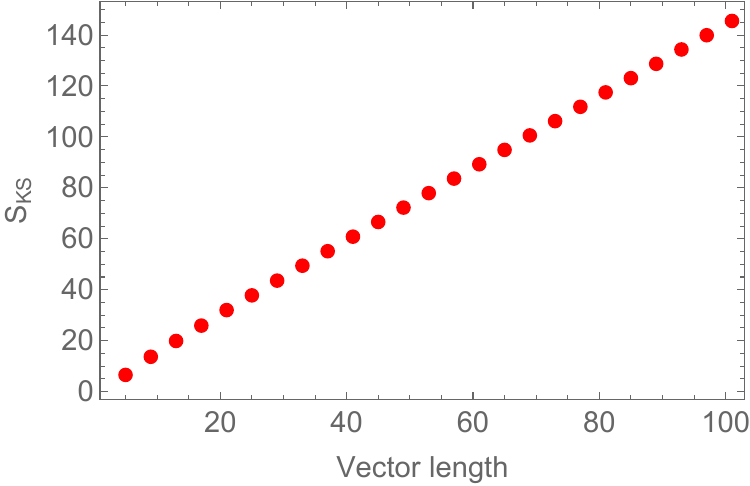}
    \caption{Kolmogorov-Sinai entropy}
    \label{fig2_subfig3}
  \end{subfigure}

  \caption{The sorted positive Lyapunov exponents for the system corresponding to the vector $h$ for the largest length $n=101$, are shown in (a), noting similar characteristics as in the case of $g$, with the more familiar drop in Lyapunov exponent magnitudes now clearly evident. As expected the K-S entropy is linear as seen in (b).}
  \label{figure_2}
\end{figure}

Moving to the second case of study, where $n=3, \dots, 101$, we must employ the connection set $S$ of a circulant graph. Using the connection set, two nodes $i$ and $j$ are connected by an edge if $(i-j)mod\, n \in S$. We call the circulant graph with connection set $S_r$, a stride $r$ graph, with the set formally defined as
\begin{equation}\label{connection_set}
    S_r = \left\{k\in \{1, 2, \dots, n-1\}\ |\ (k-1)mod\, r=0\ \text{or}\ (n-k-1)mod\, r=0 \right\}
\end{equation}
and the graph including self-loops if the connection set is of the form $\{0\} \cup S_r$. As such, in order to study how an increase in the number of nodes affects the systems, we take the example of a stride 2 circulant graph for $n=3, \dots, 101$. The stride $r$ construction is chosen as a means to increase the length of the generating vector while making sure it maintains all necessary properties. We call the circulant graph constructed in such a manner a \textit{stride 2} circulant graph, as each node is connected to the second closest node in the configuration space. The plots of the Lyapunov spectra of the system associated to the largest vector of dimension $n=101$, as well as the K-S entropy of all the systems $n=1,3,\dots,101$ can be seen in Figure~\ref{figure_2}. We omit the graphs for visual clarity. We observe the entropy increasing linearly with the system size, which serves as a numerical check of the consistency of the construction.

Clearly other constructions for increasing the vector length could be chosen and be viable. For instance, since the matrix $C$ can be viewed as the adjacency matrix of a circulant graph, inspiration from studies such as \cite{Ding2017} can be taken, and the size of the matrix $C$, and consequently the evolution matrix $M$, can be increased by turning from a fully circulant matrix to a block circulant matrix. To avoid confusion, we note that in this case the term block circulant refers to a matrix composed of blocks whose inner structure is circulant. Let $C$ and $C'$ be two circulant symmetric matrices with binary elements, which for simplicity's sake are assumed to be square of size $n\times n$. A block circulant matrix of size $2n \times 2n$ can then be constructed with $C$ and $C'$ as diagonal blocks:
\begin{equation}\label{block_circ_C}
    \tilde{C} = \begin{pmatrix}
        C & 0_{n \times n}\\
        0_{n \times n} & C'
    \end{pmatrix}
\end{equation}
Interpreting the block diagonal matrix $\tilde{C}$ as an adjacency matrix of a graph, results in two disconnected circulant graphs, where clearly a fictitious particle traveling along the nodes of one graph would be incapable of crossing over to the other. Intuitively, we would expect such a construction to fulfill all the required criteria for being a chaotic system, as constructed in \cite{Axenides2023}, since it corresponds to two separate systems. This can be proven easily, after some simple matrix multiplications, resulting in the symplectic evolution matrix
\begin{equation}\label{block_circ_evol}
    \tilde{M} = \begin{pmatrix}
        I & 0 & C & 0 \\
        0 & I & 0 & C' \\
        C & 0 & I + C^2 & 0 \\
        0 & C' & 0 & I + C'^2
    \end{pmatrix} = \begin{pmatrix}
        I_{2n \times 2n} & \tilde{C} \\
        \tilde{C} & I_{2n \times 2n} + \tilde{C}^2
    \end{pmatrix}
\end{equation}
where $I$ and $0$ are the $n$-dimensional unit and zero matrices respectively. The evolution matrix $\tilde{M}$ for the specific construction of \eqref{symp_matrix_n}, is only symplectic if the off-diagonal blocks of \eqref{block_circ_C} are zero.

\begin{figure}[htbp]
    \centering
    \includegraphics[width=0.7\linewidth]{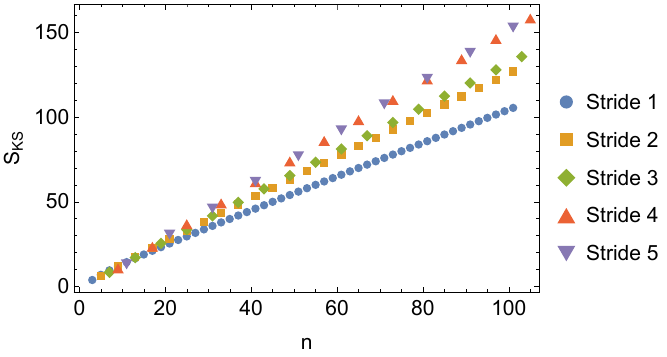}
    \caption{The K-S entropy $S_{KS}$ for five different types of graph connectivities without self loops, as a function of the length of the generating vector. It is evident that, for large $n$, the K-S entropy scaling decreases as the connectivity increases.}
    \label{figure_Sks_connectiv_1}
\end{figure}

\begin{figure}[htbp]
    \centering
    \includegraphics[width=0.7\linewidth]{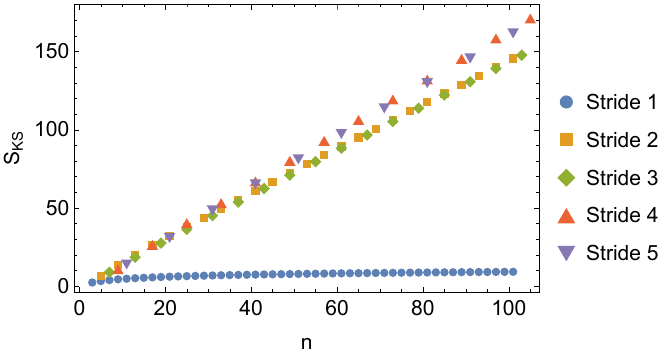}
    \caption{The K-S entropy $S_{KS}$ for five different types of graph connectivities with self loops, as a function of the length of the generating vector. A dramatic reduction in the increase of entropy production, as $n$ increases, for the fully connected case, is observed.}
    \label{figure_Sks_connectiv_2}
\end{figure}

Having studied the behaviour of the system as the number of nodes increases, we turn to the more interesting case of how the system scales with an increase to the graph connectivity, while the system size (number of nodes) remains constant. In order to do so, we build vectors following the above construction of a stride 2 circulant graph, but we now also build the stride 1, 3, 4, and 5 circulant graphs. Clearly, stride 1 corresponds to the fully connected case. We plot their K-S entropies, as well as the previous one of the step 2 case, on the same plot in Fiqures~\ref{figure_Sks_connectiv_1} and ~\ref{figure_Sks_connectiv_2} to note their different scaling.

Naively it could be expected for the entropy to increase as the connectivity of the system rises, given that a particle moving between nodes would have more possible sites to move to. Yet, as seen in Figure~\ref{figure_Sks_connectiv_1}, the opposite happens, where a decrease in connectivity leads to an increase in the entropy production of the system, with a fully connected system being characterised by the slowest increase in entropy generation. This result comes from the fact that the model's evolution is described by propagating waves, leading to cancellation of terms as a result of destructive interference. Specifically, in the stride 1 case, that is for a complete graph $K_n$, the adjacency matrix can be written as $C = \sum_{l=1}^{n-1} P^l$ and as a result it is straight forward to determine its eigenvalue spectrum $\{d_j\} = \{n-1, -1, -1, \dots, -1\}$. Consequently, the K-S entropy is given by \eqref{ks_entropy} as
\begin{equation}\label{KS_ones_no_selfloops}
    S_{KS} = \lambda_+(n-1) + (n-1) \lambda_+(-1)
\end{equation}
where the linear scaling of the entropy with $n$ is evident, alas with a shallower slope due to the degenerate -1 modes as seen in Figure~\ref{figure_Sks_connectiv_1}.

The above described destructive interference becomes even more evident if self-loops are added to the model, as seen in Figure~\ref{figure_Sks_connectiv_2}, where a fully connected graph sees almost no increase in entropy as it grows in size. In this trivial case the only non-zero Lyapunov exponent is $\lambda_{+, 0} = \ln\left(\frac{2+n^2}{2} + \frac{n}{2} \sqrt{n^2 + 4} \right)$. As a result the K-S entropy will also equal the only non-zero Lyapunov exponent of each $n$-sized system $S_{KS} = \lambda_{+, 0}$, from which it is apparent that the entropy will scale as $O(\log n)$ as $n$ grows. Interestingly, where without self loops each of systems was characterised by a different linear increase in entropy production, adding self loops leads to the cases of strides 2 and 3, as well as strides 3 and 4, to overlap.

Finally, if the reduction in entropy production is indeed caused by interference, a notable difference in $S_{KS}$ between highly structured graphs and random ones would be expected. As shown in Appendix~\ref{Appendix_b}, for a fixed node degree, graphs with internally periodic generating vectors exhibit notably lower entropy production than those with random generating vectors. This occurs because the random connectivities break the perfect interference patterns of the deterministic structures.

\subsection{Spectra of periods}
A fictitious particle's trajectory, moving on the double toroidal phase space under the action of an ACM, is always periodic, given the finite number of points on the lattice as shown by Dyson and Falk \cite{Dyson1992}, if its initial conditions are described by rational numbers. The period of the trajectory is related to the common denominator $N$ of the rational numbers and it is equal to twice the period of the Fibonacci sequence modN. In the case of many ACMs, we can determine the period of the trajectory in the multi-dimensional toroidal phase space using similar arguments as in \cite{Dyson1992}. On the other hand, initial conditions where all coordinates are rational numbers lead to non-periodic trajectories, which fill up the whole toroidal phase space. Below we describe a method to determine the periods of the periodic trajectories.

Firstly, the definition in \eqref{symp_matrix_n} can, through simple induction, be generalised to the $m$-th power, which takes the form
\begin{equation}\label{evol_matr_power}
    M^m = L^{2m} = \begin{pmatrix}
        C_{2m-1} & C_{2m} \\
        C_{2m} & C_{2m-1}
    \end{pmatrix}
\end{equation}
where $C_0 = 0_{n \times n}$, $C_1 = I_{n\times n}$, and:
\begin{equation}\label{C_matr_polynom}
    C_{m+1} = C C_m + C_{m-1} \ , \quad m=1,2,3,\dots
\end{equation}

The period of the matrix $M modN$ is the smallest integer, $T(N)$, such that 
\begin{equation}\label{evol_matrix_period_def}
    M^{T(N)} = I_{2n \times 2n}\ modN
\end{equation}
which in turn implies that
\begin{equation}
    M^{T(N)} = \begin{pmatrix}
        C_{2T(N)-1} & C_{2T(N)} \\
        C_{2T(N)} & C_{2T(N)-1}
    \end{pmatrix} modN = \begin{pmatrix}
        I_{n \times n} & 0_{n \times n} \\
        0_{n \times n} & I_{n \times n}
    \end{pmatrix} modN
\end{equation}
leading to the relations
\begin{equation}\label{periodic_conditions}
    \begin{split}
        C_{2T(N)-1} &= I_{n \times n}\ modN \\
        C_{2T(N)} & = 0_{n \times n}\ modN
    \end{split}
\end{equation}
which must both hold simultaneously.

We can show that the recursion relation \eqref{C_matr_polynom} is similar to the one defining the Fibonacci polynomials \cite{Philippou2002},
\begin{equation}\label{fib_poynomials_def}
    F_{m+1}(x) = x F_m(x) + F_{m-1}(x)
\end{equation}
with $x \in \mathbb{R}$, and $F_0(x) = 0$ and $F_1(x) = 1$. Consequently, we can make use of the known formula of the Fibonacci polynomials' form, to determine the general form of the matrices $C_m$
\begin{equation}
    C_m \equiv F_m(C) = \sum_{j=0}^{\left[ \frac{m-1}{2}\right]} \begin{pmatrix}
        m - j - 1 \\
        j
    \end{pmatrix} C^{m - 2j -1}
\end{equation}
with $[.]$ denoting the integer part. In the above relation we have replaced the argument $x$ of the Fibonacci polynomials by the matrix $C$.

We now present some numerical results to better visualise how the period of the evolution matrix $M$ scales both with the resolution of the phase space, through the small pieces of size $1/N$, as well as with the number $n$ of the number of ACMs.

\begin{figure}[htbp]
    \centering
    \includegraphics[width=0.9\textwidth]{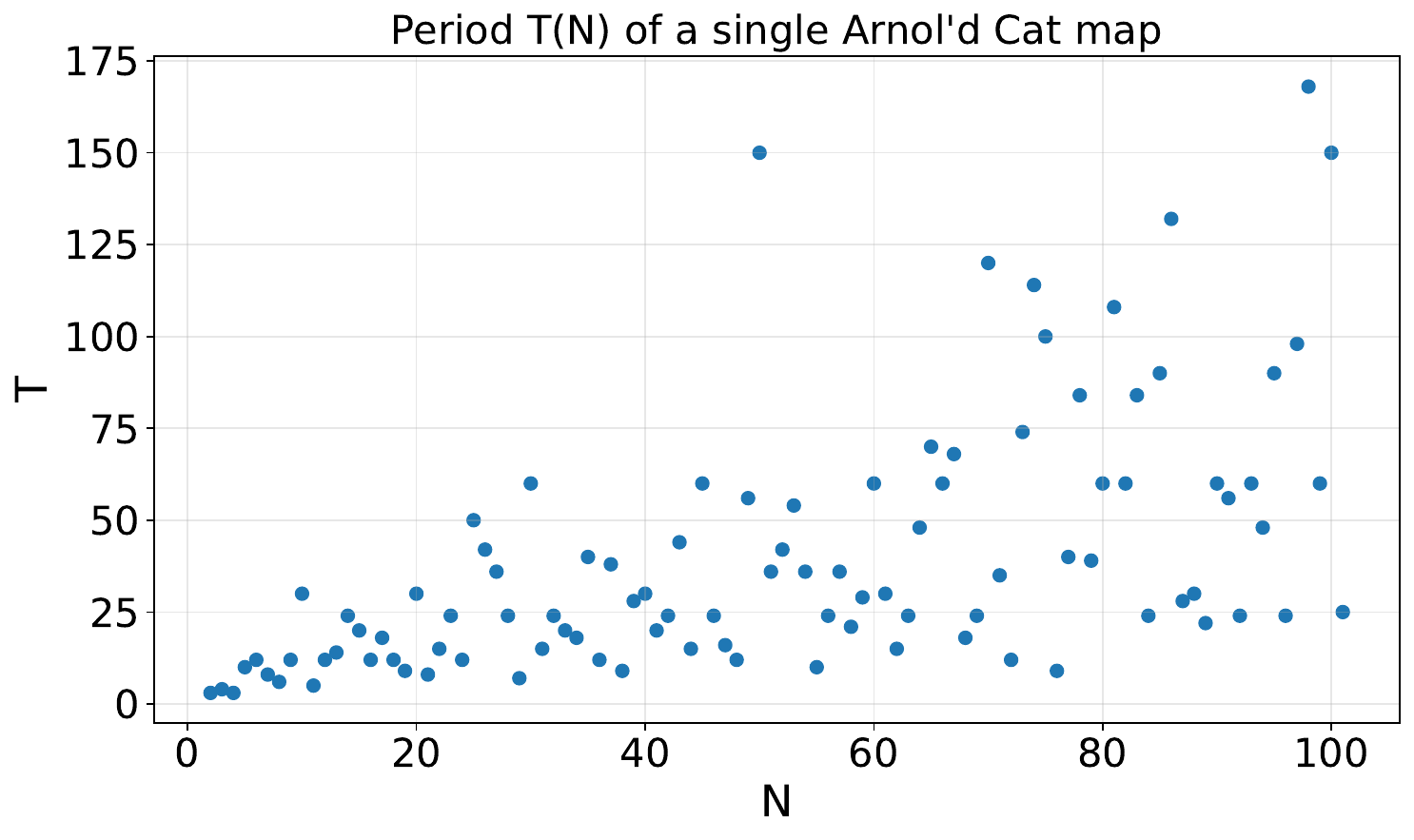}
    \caption{The period spectra, as function of the phase space resolution $T=T(N)$, of the single ACM.}
    \label{figure_single_ACM}
\end{figure}

\begin{figure}[htbp]
    \centering
    \includegraphics[width=0.9\textwidth]{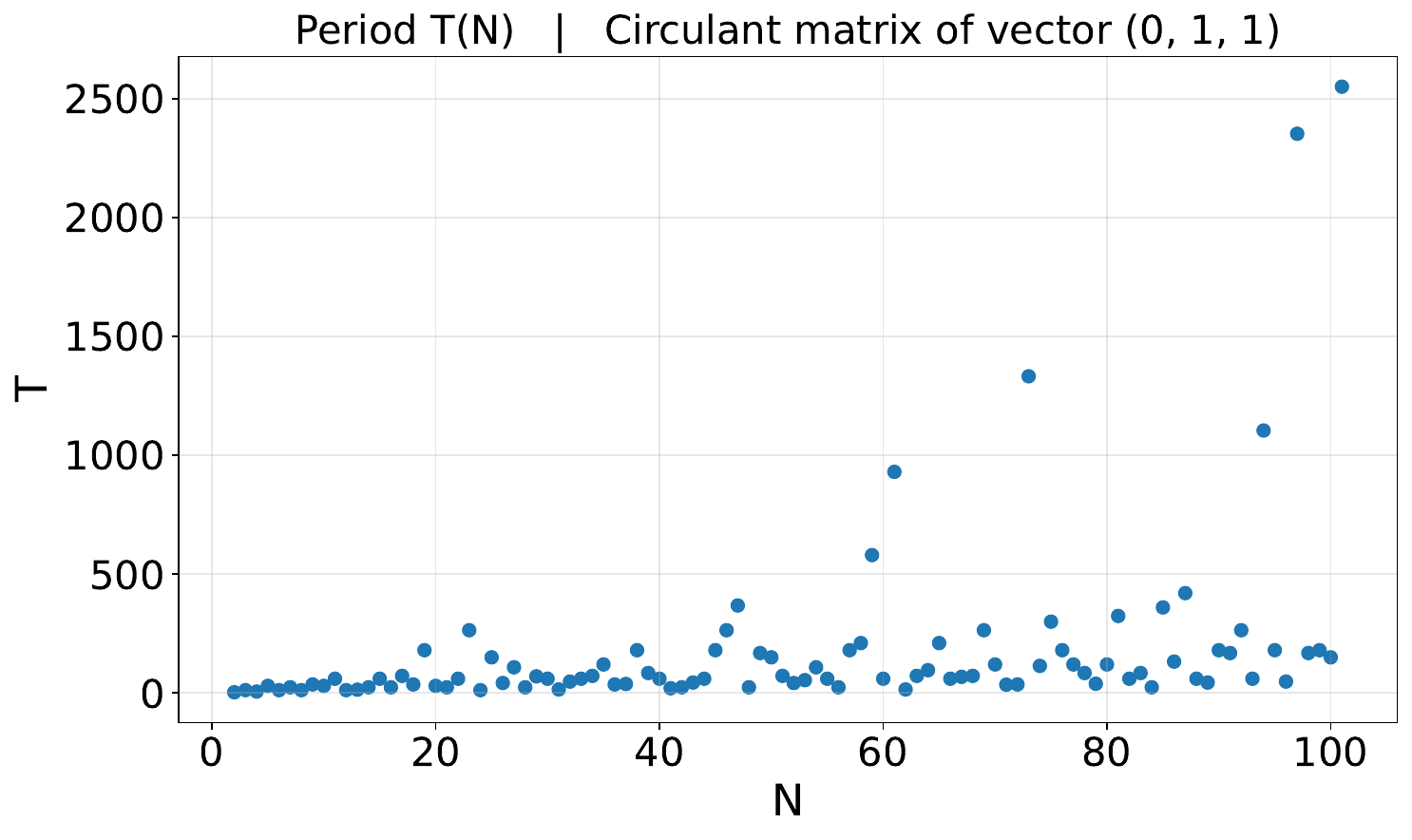}
    \caption{The period spectra, as function of the phase space resolution $T=T(N)$, of the evolution matrix $M$ for three coupled ACMs. The matrix $M$ is constructed through the circulant matrix $C$ which itself is generated by the vector $(0, 1, 1)$.}
    \label{figure_three_ACMs}
\end{figure}

\begin{figure}[htbp]
    \centering
    \includegraphics[width=0.9\textwidth]{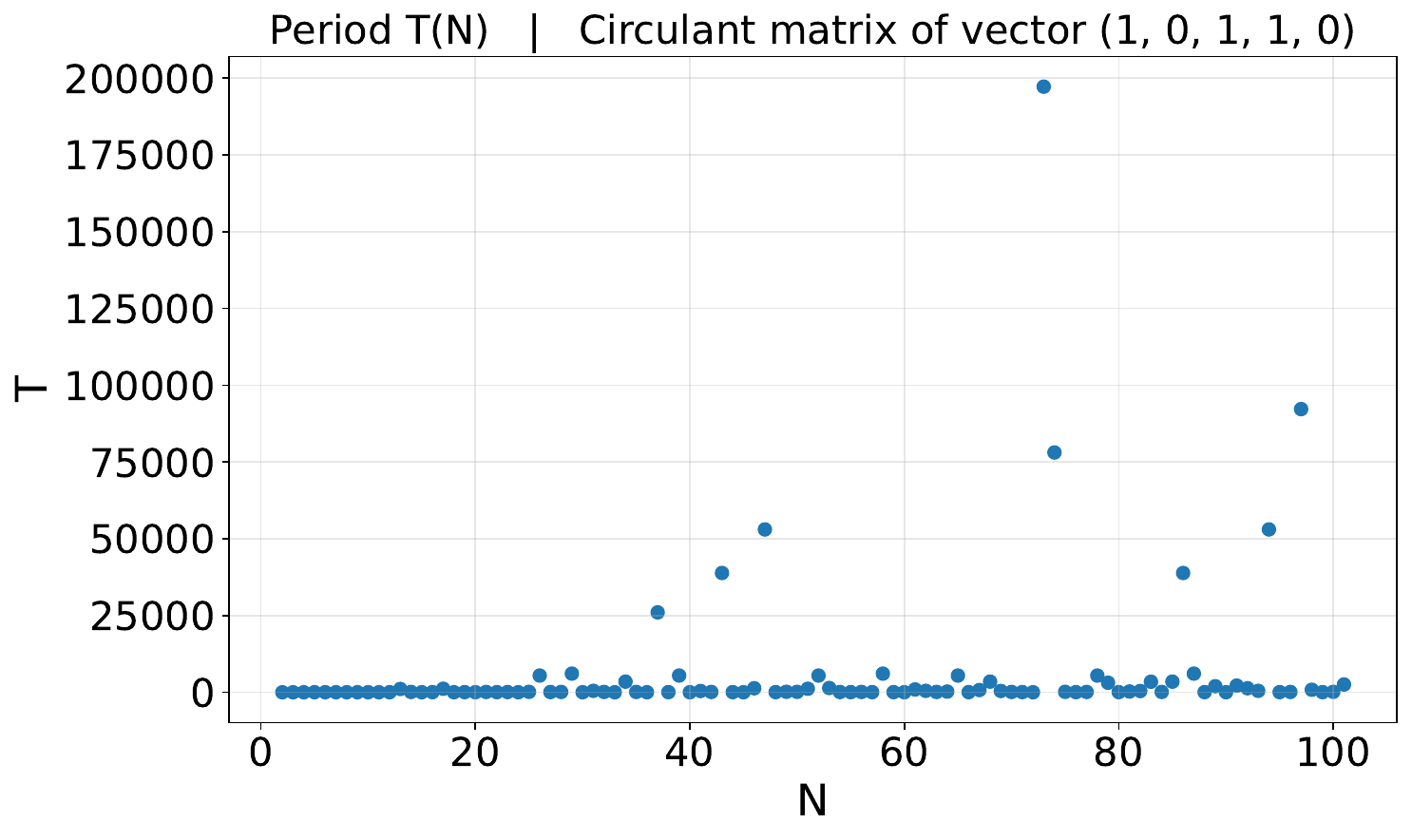}
    \caption{The period spectra, as function of the phase space resolution $T=T(N)$, for an evolution matrix $M$ of size 10. The matrix $C$ is generated by the vector $(1, 0, 1, 1, 0)$.}
    \label{figure_five_ACMs}
\end{figure}

\begin{figure}[htbp]
    \centering
    \includegraphics[width=0.9\textwidth]{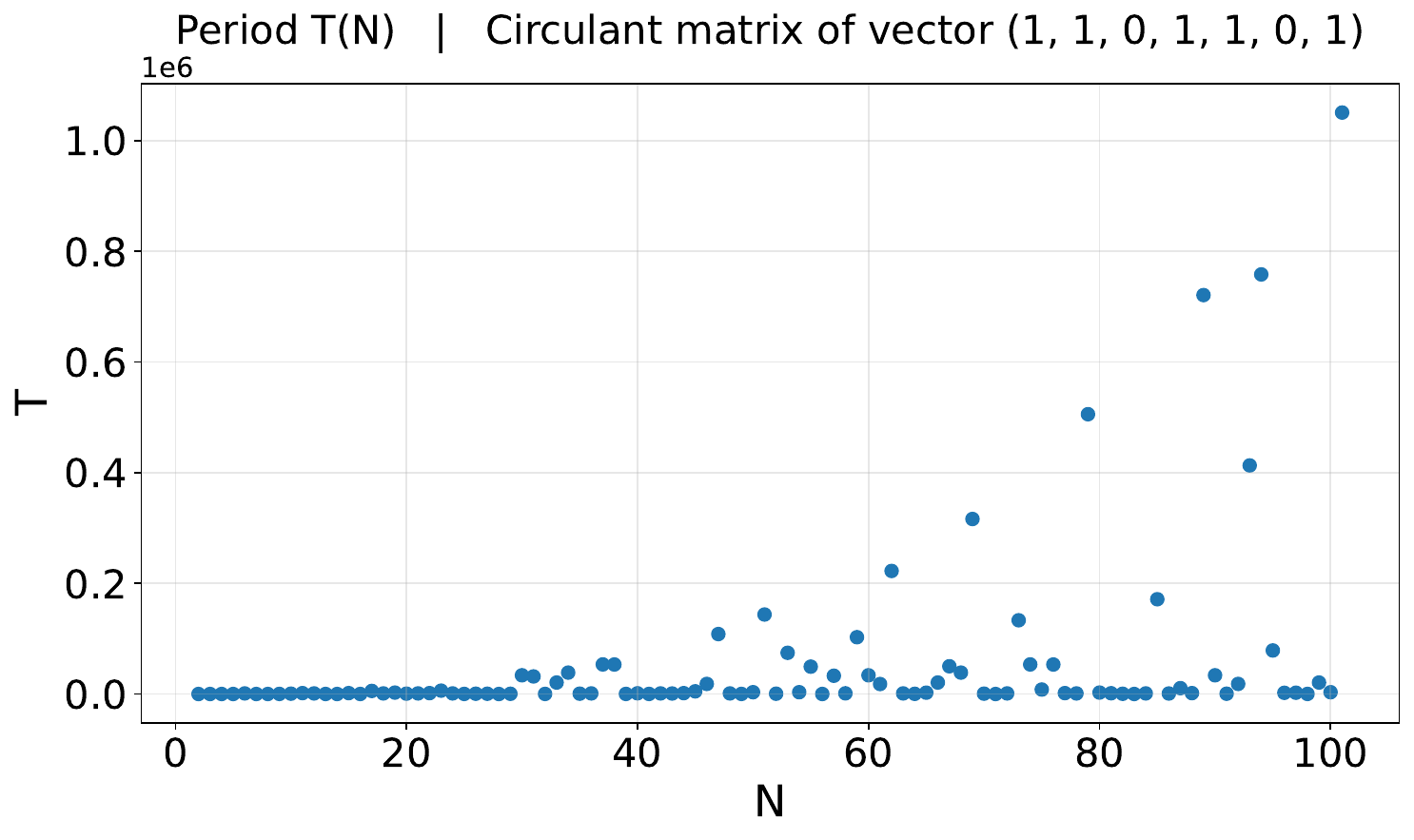}
    \caption{The period spectra, as function of the phase space resolution $T=T(N)$, for an evolution matrix $M$ of size 14. The matrix $C$ is generated by the vector $(1, 1, 0, 1, 1, 0, 1)$.}
    \label{figure_seven_ACMs}
\end{figure}

First the period $T(N)$ as a function of the torus resolution $N$ is shown in Figures ~\ref{figure_single_ACM}, ~\ref{figure_three_ACMs}, ~\ref{figure_five_ACMs}, and ~\ref{figure_seven_ACMs} for the case of a singular ACM, as well as for matrices $C$ of dimension $3,5$ and $7$.

\begin{figure}[htbp]
  \centering
    \centering
    \includegraphics[width=\textwidth]{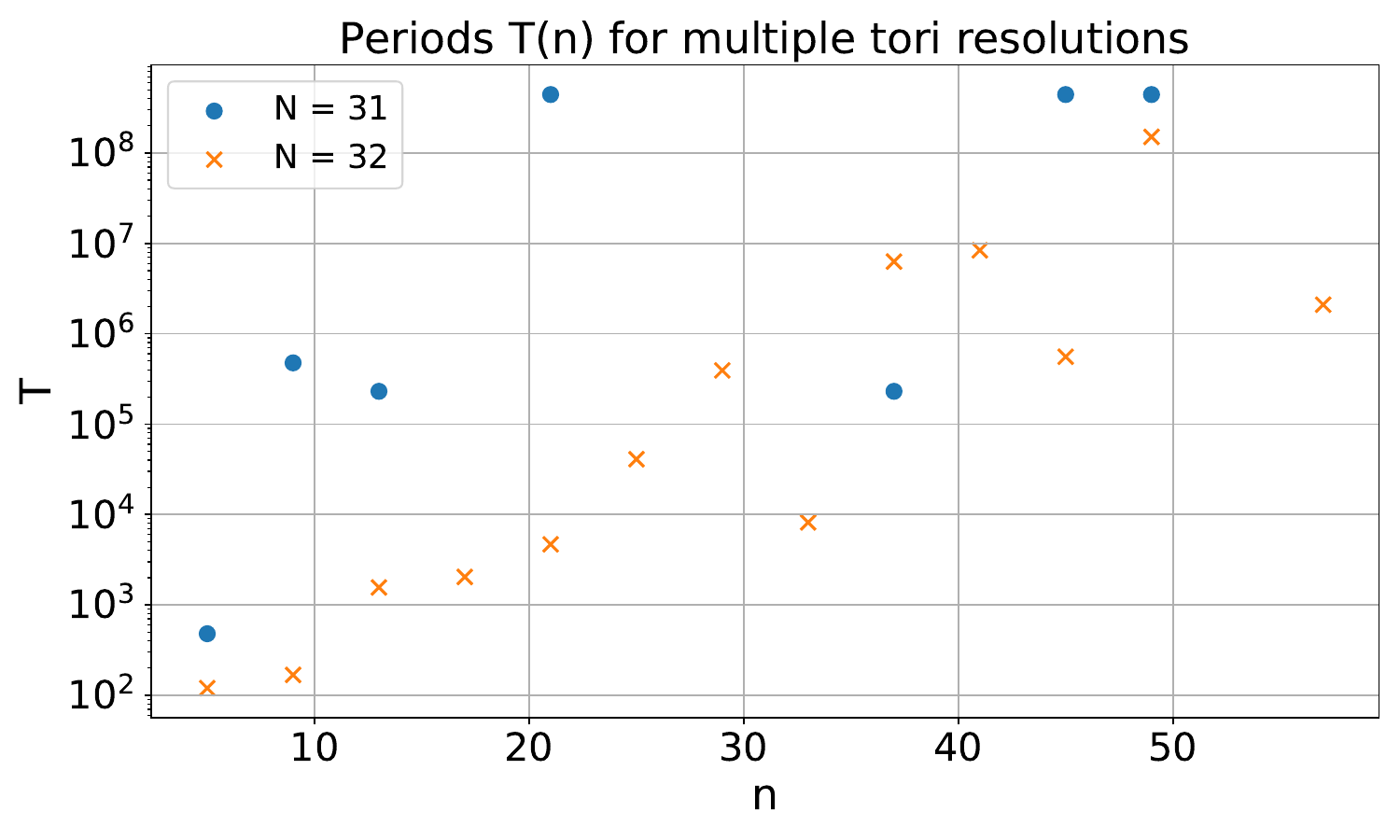}
  \caption{A plot of the period spectra $T(n)$ as a function of the number of coupled ACMs, for $N=\{31,32\}$ is shown. The matrix $C$ for these two systems is generated by a vector with periodic structure. Specifically, the vector is initialised as $(0,0,1)$ and then the sequence $(0,1)$ is repeatedly added to increase its size up to a maximum length of $50$. At each step the vector grows by $4$, as we also add the reverse sequence $(1,0)$ to the latter half of the vector to observe the symmetry constraints. Due to computational constraints, an upper limit of $T=10^9$ for the period was imposed, and data points exceeding this limit are excluded (resulting in missing $T(n)$ values for $N=31$). A vector generated by random bits, still observing the needed symmetry constraints, was also explored but no significant difference in behaviour was observed.}
  \label{figure_5}
\end{figure}

Finally, we provide a plot of the period $T(n)$ as a function of the number of the degrees of freedom in Figure \ref{figure_5} for $N=31$ and $N=32$. Other values of $N$ were tested, such as non-prime odds and non-power of two evens, and similar behaviour was observed.

For values of $N \neq 2^s, \ s=1,2 \dots$ the period spectra retains its chaotic behaviour, where even for $n \simeq 15$ a single addition of a node can change the period by many orders of magnitude. On the other hand for $N = 2^s,\ s=1,2,\dots$, the period spectra, although still random in nature, is greatly constrained in comparison. Notably, in this case the system behaves markedly different for even and odd sizes of the matrix $C$. 

In the even case the periods do not depend on the number of coupled ACMs, that is the size of the matrix $C$, but behave similarly to the case of the single ACM exhibiting the same upper and lower bounds:
\begin{equation}\label{period_bounds_power_two}
    \frac{3N}{4} \leq T \leq 3N
\end{equation}
Specifically the periods of the system are dictated by $C^2$ (as will be shown bellow), which as the square of the adjacency matrix of a graph, corresponds to neighbours of only distance 2 communicating with each other leading to the circulant graph of coupled ACMs being split in two disconnected graphs of nodes on even and odd sites. Given the block form of the evolution matrix $M$, it can be shown with simple algebra calculations that the eigenvalues obey the equation
\begin{equation}\label{M_C_eigenva}
    \lambda + \frac{1}{\lambda} = 2 + \mu^2
\end{equation}
where $\lambda$ the eigenvalues of $M$ and $\mu$ the eigenvalues of $C$, which taking $modN$ gives
\begin{equation}\label{M_C_eigenva_mod2}
    \lambda + \frac{1}{\lambda} = \mu^2
\end{equation}
making evident the dependence of the period on $\mu^2$, since finding the period $T$ such that $M^T = ImodN$ simply corresponds to finding $T$ such that all eigenvalues $\lambda^T = 1modN$. Also, although we cannot provide an analytic form for an arbitrary period, in this case if the period for a specific number $n$ of coupled ACMs for $N=2$ is known, then the period for the same $n$ but larger $N$ can be determined by
\begin{equation}\label{period_power2_form}
    T(N;n) = \frac{N}{2} T(2;n) \ ,\quad N=2^s, \ s=2,3,4,\dots
\end{equation}
for a set number $n$ of ACMs.

In the odd case the graph does not split in two and as such maximal spatial mixing is achieved. An analytic form for the base period could potentially be determined, given the similarities to Cellular Automaton, but that is beyond the scope of this work. The upper bound will scale exponentially with the number $n$ of ACMs, and the lower bound will scale linearly with $n$ as the information still needs to propagate the entire length of the graph given the maximal spatial mixing.

\section{Conclusion and future directions}\label{section_conclusion}
In this work, we have generalised the model of symplectically interacting Arnol'd cat maps (ACMs), of reference \cite{Axenides2023}, into a model that describes the chaotic properties of general circulant graphs, by placing an ACM on every node of the graph and coupling them based on the connectivity of the graph. By demanding that the system's evolution matrix is symplectic, the coupling matrix naturally becomes symmetric, allowing it to be interpreted directly as the adjacency matrix of a circulant undirected graph. This method offers a new approach to solvable models of chaotic dynamical systems where the discrete Fourier transform enables exact analysis. 

A novel observation of this study is the counterintuitive property of the reduction of the KS entropy with the increase in graph connectivity, due to the circulant nature of the graph. Furthermore, our numerical simulations for the periods of this system reveals that, when the spatial resolution of the tori $N$ is a power of two, coupling an even number of ACMs leads to period spectra similar to that of a single ACM.

There are two natural avenues for future research, depending on whether or not the circulant structure of the matrices is sought to be maintained. If the aim is to generalise to non-circulant matrices, then translation invariance, in our framework, needs to be abandoned and the solvability of the model is lost. In this case the coupling matrix of \eqref{near_neighb_coupl_matr}, for an arbitrary range of interactions, takes the form:
\begin{equation}\label{gen_coupl_matr_non_trans}
    C = G_0 I_{n \times n} + \sum_{l=1}^{[n/2]}\left( P^l G_l + G_l \left( P^T \right) \right)
\end{equation}

Two other directions that can be taken are the generalisations to higher dimensional surfaces, or to time-varying connectivities. Generalisation to higher dimensional surfaces can be achieved by extending the coupling matrix to a Block Circulant matrix with Circulant Blocks (BCCB). Because BCCB matrices maintain the necessary symmetry to ensure a symplectic evolution matrix (as shown in Appendix \ref{ProofBCCB}) and are diagonalized by the multidimensional Discrete Fourier Transform, they provide an analytically exact method to investigate whether the destructive interference that suppresses entropy generation in 1D networks persists or scales differently across two or three spatial dimensions. On the other hand, generalising to time-varying connectivities is straight forward to accomplish, by taking the product of multiple symplectic evolution matrices, which will also be symplectic. As such is would be worthwhile to explore whether dynamically switching between circulant topologies of differing connectivities induces resonances or further limits entropy production compared to a static topology.

Additionally, the number-theoretic properties governing the period spectra on these finite toroidal spaces warrant deeper exploration. Our numerical results indicated that when the size of the phase space $N=2^s$, the allowable periods are heavily restricted, even for an odd number of nodes. This relates to the algebraic structure of the symplectic group, $Sp_{2n}(\mathbb{Z}_{2^s})$, which could be explored using number theoretic methods established for multidimensional toral maps \cite{Kurlberg2000}. While it is already established that the period growth for a single ACM follows set lower and upper bounds, some of which pass on to the case of $N=2^s, \ s=1,2,\dots$ for multiple coupled ACMs, future work could seek to derive rigorous analytical bounds for the period $T(N)$ for any $N$.

Finally, due to the symplecticity of the evolution matrix of the coupled ACMs it is possible to study the quantum mechanical chaotic properties of this model. Given the model's inherent connections to translation operators on the toroidal lattice in Finite Quantum Mechanics, quantizing these coupled maps could provide valuable insights into many-particle quantum chaos. Specifically, future studies could investigate whether the specific graph topologies that constrain classical entropy production similarly limit the growth of quantum entanglement entropy, or if these symmetric circulant structures can be leveraged to facilitate perfect quantum state transfer \cite{BASIC2009}.

\section*{Acknowledgments}
The authors wish to thank Dr. Ioannis Tsohantjis for their insightful conversations and valuable feedback during the development of this work.

\appendix
\section{Proof of Symplecticity for 2D BCCB Evolution Matrix}\label{ProofBCCB}
Let the symmetric block circulant matrix with circulant blocks (BCCB) be defined as:
\begin{equation}
\tilde{C}_{2n \times 2n} = \begin{pmatrix} C_{n \times n} & B_{n \times n} \\ B^T_{n \times n} & C_{n \times n}' \end{pmatrix}
\end{equation}
where $C_{n \times n} = C_{n \times n}^T$ and $C_{n \times n}' = (C_{n \times n}')^T$.
The evolution matrix $M$, of dimension $4n \times 4n$, is constructed as the square of the anti-symplectic Fibonacci matrix, which yields:
\begin{equation}
M = \begin{pmatrix} I & \tilde{C} \\ \tilde{C} & I + \tilde{C}^2 \end{pmatrix} 
= \begin{pmatrix} 
I & 0 & C & B \\ 
0 & I & B^T & C' \\ 
C & B & I + C^2 + BB^T & CB + BC' \\ 
B^T & C' & B^TC + C'B^T & I + B^TB + C'^2 
\end{pmatrix}
\end{equation}
The symplectic form $J$ in this $4 \times 4$ block space is:
\begin{equation}
J = \begin{pmatrix} 
0 & 0 & -I & 0 \\ 
0 & 0 & 0 & -I \\ 
I & 0 & 0 & 0 \\ 
0 & I & 0 & 0 
\end{pmatrix}
\end{equation}
For notational simplicity, we have ignored the explicit annotation of the dimension of each matrix, since all matrices $I, B, C, C'$ are of dimension ${n \times n}$.

In order to prove symplecticity, we calculate first
\begin{equation}
JM = \begin{pmatrix} 
-C & -B & -(I + C^2 + BB^T) & -(CB + BC') \\ 
-B^T & -C' & -(B^TC + C'B^T) & -(I + B^TB + C'^2) \\ 
I & 0 & C & B \\ 
0 & I & B^T & C' 
\end{pmatrix}
\end{equation}
and then with some tedious, yet simple, multiplications it can be shown that:
\begin{equation}
M^T J M = M J M = \begin{pmatrix} 0 & 0 & -I & 0 \\ 0 & 0 & 0 & -I \\ I & 0 & 0 & 0 \\ 0 & I & 0 & 0 \end{pmatrix} = J
\end{equation}

As such, so long as the structure of the coupling matrix is BCCB, the two-dimensional lattice of nodes also leads to symplectic evolution. Geometrically this construction corresponds to multiple circulant graphs with 1D ring topologies, stacked on top of each other and connected as dictated by the off-diagonal blocks of the coupling matrix $\tilde{C}$. Enforcing translation invariance on the second dimension too, leads to a topology of a 2D torus.

\section{Entropy dependence on structural periodicity}\label{Appendix_b}
\begin{figure}[h]
    \centering
    \includegraphics[width=0.9\linewidth]{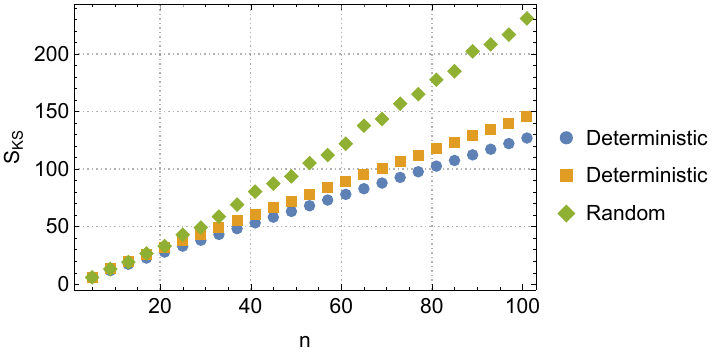}
    \caption{Kolmogorov-Sinai entropy ($S_{KS}$) as a function of the number of nodes ($n$) for circulant graphs of the same connection density (stride 2). The plot compares deterministically structured connectivities against random symmetric connectivities of the identical node degree. Specifically, one of the graphs labeled deterministic, has the connection set $S_2$, as defined in \eqref{connection_set}, whereas the other has the connection set $S_2' = \left\{k\in \{1, 2, \dots, n-1\}\ |\ k\,mod\, r=0\ \text{or}\ (n-k)mod\, r=0 \right\}$. Data points for the random configurations represent the mean average over ten runs to account for statistical variance. The random configurations consistently exhibit higher entropy production, demonstrating that the strict structural periodicity of the deterministic graphs induces spectral degeneracies that suppress chaotic expansion. Similar behaviour was observed for all different node degrees tested.}
    \label{figure_appendix}
\end{figure}

\clearpage
\bibliographystyle{unsrturl}
\bibliography{references}
\end{document}